\documentclass[12pt]{article}

\usepackage{a4wide}
\usepackage{amssymb}
\usepackage{amsfonts}
\usepackage{amsmath}
\input xy
\xyoption{arrow} \xyoption{matrix}

\date{}

\newtheorem{proposition}{Proposition}[section]
\newtheorem{theorem}[proposition]{Theorem}
\newtheorem{lemma}[proposition]{Lemma}

\newtheorem{corollary}[proposition]{Corollary}

\def\GK{{\rm  GK}\,}

\def\der{\partial }

\def\nFM0{{\nu }_{F,M_0}}
\def\nFN0{{\nu }_{F,N_0}}
\def\nGN0{{\nu }_{G,N_0}}

\def\N0{ {\bf N}_0 }

\def\t{\otimes}

\def\ra{\rightarrow}

\def\Xpm{X^{\pm }}

\def\s{\sigma}

\def\l1{{\lambda}_1}

\def\a{\alpha}
\def\a0{ {\alpha }_0}
\def\a1{ {\alpha }_1}

\def\l{\lambda}

\def\nFGM0{{\nu }_{F,G,M_0}}

\def\nFN0{{\nu}_{F,N_0}}

\def\sm{{\sigma}^m}

\def\sm1{{\sigma}^{-1}}

\def\smtp1{{\sigma}^{-t+1}}

\def\S1{S^{-1}}

\def\Xpm1{X^{\pm 1}_1}

\def\sPM1{{\sigma }^{\pm 1}}
\def\sMP1{{\sigma }^{\mp 1 }}

\def\d{\delta}

\def\di{{\rm d.ind}}

\def\L{\Lambda}

\def\CA{{\cal A}}

\def\CD{{\cal D}}

\def\Ytm1{Y^{t-1}}
\def\Yim1{Y^{i-1}}

\def\Aut{{\rm Aut}}

\def\ad{{\rm ad }}
\def\dim{{\rm dim }}

\def\ker{ {\rm ker } }

\def\D{ \Delta }

\def\SL2Z{ {\rm SL}_2({\bf Z}) }

\def\Gp1{ G^{1 , 1 } }
\def\P11{ P^{-1 , 1 } }
\def\Pp1{ P^{1 , 1 } }

\def\nCLsr{{}^\nu\kern-2pt {\cal L}^{\sigma , \rho  }}
\def\nP{{}^\nu \kern-2pt P}
\def\nL{{}^\nu\kern-2pt L}
\def\nLL{{}^\nu\kern-2pt \Lambda}
\def\nPsr{{}^\nu\kern-2pt P^{\sigma , \rho  }}
\def\nLsr{{}^\nu\kern-2pt L^{\sigma , \rho  }}
\def\nuCL{{}^\nu\kern-2pt  {\cal L}}
\def\nCLsr{{}^\nu\kern-2pt {\cal L}^{\sigma , \rho  }}
\def\nCL1m{{}^\nu\kern-2pt {\cal L}^{-1 , 1  }}
\def\x1nu{x^\frac{1}{\nu}}
\def\xm1nu{x^{-\frac{1}{\nu}}}

\def\ra{\rightarrow }

\def\CP{ {\cal P}}

\def\nAM0{{\nu }_{{\cal A},M_0}}
\def\nAN0{{\nu }_{{\cal A},N_0}}

\def\CP{ {\cal P }}
\def\det{ {\rm det }}
\def\ad{ {\rm ad }}

\def\gn{\mathfrak{n}}
\def\gm{\mathfrak{m}}
\def\gp{\mathfrak{p}}

\def\Spec{{\rm Spec}}

\def\di!{\frac{\der^i}{i!}}
\def\dik!{\frac{\der^k_i}{k!}}
\begin{document}

\author{V. V. \  Bavula 
}

\title{The inversion formula for automorphisms of the Weyl algebras and  polynomial algebras}

\maketitle
\begin{abstract}
Let $A_n$ be the $n$'th {\em Weyl algebra} and $P_m$ be a
polynomial algebra in $m$ variables over a field $K$ of
characteristic zero. The following characterization of the
algebras $\{ A_n\t P_m\}$ is proved: {\em an algebra $A$ admits a
 finite set $\d_1, \ldots , \d_s$ of commuting locally nilpotent derivations with generic
 kernels and $\cap_{i=1}^s\ker (\d_i)=K$ iff $A\simeq A_n\t P_m$
 for some $n$ and $m$ with $2n+m=s$, and vice versa}. {\em The inversion
 formula} 
 for  automorphisms of the algebra   $A_n\t P_m$ (and for $\widehat{P}_m:= K[[x_1, \ldots , x_m]]$) has
 found (giving a {\em new} inversion formula even for
 polynomials). Recall   that (see \cite{BCW}) {\em given $\s \in \Aut_K(P_m)$, then}  $\deg \,
 \s^{-1}\leq (\deg \, \s )^{m-1}$ (the proof is {\em algebro-geometric}).
We extend this result (using [non-holonomic] $\CD$-{\em modules}):
{\em given $\s \in \Aut_K(A_n\t P_m)$, then }$\deg \,
 \s^{-1}\leq (\deg \, \s )^{2n+m-1}$. Any automorphism  $\s \in
 \Aut_K(P_m)$ is determined by its {\em face} polynomials \cite{McKayWang88Inv}, a
 {\em similar} result is proved for $\s \in \Aut_K(A_n\t P_m)$.

One can amalgamate two old open problems ({\bf the Jacobian
Conjecture} and {\bf the Dixmier Problem}, see \cite{Dix} problem
1) into a single question, ({\bf JD}): {\em is a $K$-algebra
endomorphism $\s :A_n\t P_m\ra A_n\t P_m$  an algebra automorphism
provided $\det (\frac{\der \s (x_i)}{\der x_j})\in
K^*:=K\backslash \{ 0\}$?} $(P_m=K[x_1, \ldots , x_m])$. It
follows immediately from the inversion formula that {\em this
question has an affirmative answer iff both conjectures have} (see
below) [{\em iff one of the conjectures has a positive answer} (as
it  follows from the recent paper \cite{Bel-Kon05JCDP})].

 {\em Mathematics subject classification
2000: 13N10, 13N15, 14R15, 14H37,  16S32.}
\end{abstract}


\section{Introduction}

The following notation will remain {\bf fixed} throughout the
paper (if it is not stated otherwise): $K$ is a field of
characteristic zero (not necessarily algebraically closed), module
means a {\em left} module, 
$A_n=\oplus_{\alpha \in \mathbb{N}^{2n}}Kx^\alpha$ is the $n$'th
{\em Weyl algebra} over $K$, $P_m=\oplus_{\alpha \in
\mathbb{N}^{m}}Kx^\alpha$ is a {\em polynomial} algebra over $K$,
$A:= A_n\t P_m= \oplus_{\alpha \in \mathbb{N}^{s}}Kx^\alpha$,
$s:=2n+m$, is the Weyl algebra with polynomial coefficients where
$x_1, \ldots, x_s$ are the canonical generators  for $A$ (see
below). Any $K$-algebra automorphism $\s \in \Aut_K(A)$ is
uniquely determined by the elements $x_i':= \s (x_i)=\sum_{\alpha
\in \mathbb{N}^s}\l_\alpha x^\alpha$, $i=1, \ldots , s$,
$\l_\alpha \in K$, and so does its inverse,
$\s^{-1}(x_i)=\sum_{\alpha \in \mathbb{N}^s}\l_\alpha' x^\alpha$,
$i=1, \ldots , s$.

{\bf What is \underline{the} inversion formula for $\s \in
\Aut_K(A_n\t P_m)$?} A natural (shortest)   answer to this
question is a formula for the coefficients
$\l_\alpha'=\l_\alpha'(\l_\beta )$ {\em like the inversion formula
(the Kramer's formula) in the linear polynomial case}: given
$x'=Ax$ where $A=(a_{ij})\in {\rm GL}_m(K)$ (i.e.
$x_i'=\sum_{j=1}^ma_{ij}x_j$ where $a_{ij}=\frac{\der x_i'}{\der
x_j}$) then
$$x=A^{-1}x'=(\frac{\der x_i'}{\der x_j})^{-1}x'= (\det \, A)^{-1}
(\D_{ij}) x'$$ where $\D_{ij}$ are complementary minors for the
matrix $(\frac{\der x_i'}{\der x_j})$, they are linear
combinations of products of partial derivatives $\frac{\der
x_i'}{\der x_j}$. So, the inversion formula in general situation
is  a formula, $\l_\alpha'=\l_\alpha'(\l_\beta )$, where {\em
only} additions and multiplications are allowed of `partial
derivatives' of the elements $x'$ (taking partial derivatives
`correspond' to operation of taking coefficients of $x'$). So,
{\em the} inversion formula is the most {\em economical} formula
(the point I want to make is that $x=\frac{1}{2}x'$ is the
inversion formula for the equation $x'=2x$ but
$x=\frac{1}{2}(x'+2\int_0^1f(t)dt+2\dim_K{\rm
Ext}^i_B(M,N))-\int_0^1f(t)dt-\dim_K{\rm Ext}^i_B(M,N)$ is `not').

Theorem \ref{i8Nov05} gives {\em the inversion formula} for an
automorphism $\s \in \Aut_K(A_n\t P_m)$. Theorem \ref{hPi8Nov05}
gives a {\em similar formula} for an automorphism $\s \in
\Aut_K(K[[x_1, \ldots , x_m]])$. For {\em another} inversion
formula for $\s \in \Aut_K(P_m)$ see \cite{BCW},
\cite{Adj-vdEssen90}.

{\bf The degree of $\s^{-1}$  where $\s \in \Aut_K(A_n\t P_m)$}.
 We extend the following result which according to the comment made
on p. 292, \cite{BCW}: `was ``well-known'' to the classical
geometers' and `was communicated to us [H. Bass, E. H. Connell, D.
Wright] by Ofer Gabber. ... He attributes it to an unrecalled
colloquium lecture at Harward.'

\begin{theorem}\label{}
\cite{BCW}, \cite{Rus-Win84}. Given $\s \in \Aut_K(P_m)$, then
$\deg \,
 \s^{-1}\leq (\deg \, \s )^{m-1}$.
\end{theorem}
The proof of this theorem is {\em algebro-geometric} (see
\cite{Adj-Win05} for a generalization of this result for certain
varieties)). We extend this result (see Section \ref{DEGINVA}).

\begin{theorem}\label{intr9Nov05}
Given $\s \in \Aut_K(A_n\t P_m)$. Then $ \deg \, \s^{-1}\leq (\deg
\, \s )^{2n+m-1}$.
\end{theorem}
{\em Non-holonomic} $\CD$-modules are used in the proof (it looks
like this is one of the first instances where {\em non-holonomic}
$\CD$-modules are of real use).

{\bf The algebras $\{ A_n\t P_m\}$ as a class}. Theorem
\ref{Wpolchar} gives a characterization  of the algebras  $\{
A_n\t P_m\}$ as a class via commuting sets of locally nilpotent
derivations: {\em an algebra $A$ admits a
 finite set $\d_1, \ldots , \d_s$ of commuting locally nilpotent derivations with generic
 kernels and $\cap_{i=1}^s\ker (\d_i)=K$ iff $A\simeq A_n\t P_m$
 for some $n$ and $m$ with $2n+m=s$, and vice versa} (the kernels $\ker (\d_i)$
  are {\em generic} if the intersections $\{ \cap_{i\neq j}\ker (\d_i)\, | \, j=1, \ldots , s\}$
   are {\em distinct}).

{\bf Left and right faces of an automorphism $\s \in \Aut_K(A_n\t
P_m)$}.  Let $P_m=K[X_1, \ldots , X_m]$ be a polynomial algebra.
For each $i=1, \ldots , m$, the {\em algebra epimorphism}
$f_i:P_m\ra P_m/(X_i)$, $p\mapsto p+(X_i)$,  is called the {\em
face homomorphism}. J. McKay and S. S.-S. Wang
\cite{McKayWang88Inv} proved: {\em given $\s , \tau  \in
\Aut_K(P_m)$ such that $f_i\s =f_i\tau $, $i=1, \ldots , m$, then}
$\s =\tau$. So, an automorphism $\s \in \Aut_K(P_m)$ is completely
determined by its {\em faces} $\{ f_i\s \, | \, i=1, \ldots , m\}$
or equivalently by its {\em face polynomials} $\{ f_i\s (X_j)\, |
\, i, j=1, \ldots , m\}$ since each $f_i\s $ is an {\em algebra
homomorphism}.

For the algebra $A:= A_n\t P_m=K\langle x_1, \ldots , x_s\rangle$,
$s:=2n+m$, (where $x_1, \ldots , x_s$ are the canonical
generators) we have {\em left faces} $l_i: A\ra A/Ax_i$, $
a\mapsto a+Ax_i$, and {\em right  faces} $r_i: A\ra A/x_iA$, $
a\mapsto a+x_iA$, $i=1,\ldots , s$. These are homomorphisms of
{\em left} and {\em right} $A$-modules rather then homomorphisms
of algebras (if $x_i\in P_m$ then $l_i=r_i$ is an {\em algebra
homomorphism}).

Theorem \ref{27Nov05} states that: {\em given $\s , \tau \in
\Aut_K(A)$ such that  $r_i\s =r_i\tau$, $i=1, \ldots , s$ then $\s
=\tau $ (similarly, $l_i\s =l_i\tau$, $i=1, \ldots , s$, imply $\s
=\tau $)}.


\section{The Inversion Formula}
In this section, the inversion formula (Theorem \ref{i8Nov05}) is
given.

Let $A$ be an algebra over a field $K$ and let $\d $ be a
 $K$-derivation of the algebra $A$. For any elements $a,b\in A$
 and a natural number $n$, an easy induction argument gives
 $$ \d^n(ab)=\sum_{i=0}^n\, {n\choose i}\d^i(a)\d^{n-i}(b).$$
 It follows that the kernel $A^\d:=\ker \, \d $ of $\d $ is a
 subalgebra (of {\em constants} for $\d $) of $A$ and the union of the
 vector spaces $N:=N(\d ,A)=\cup_{i\geq 0}\, N_i$ is a positively
 {\em filtered} algebra ($N_iN_j\subseteq N_{i+j}$ for
 all $i,j\geq 0$).
 Clearly, $N_0= A^\d$ and $N:=\{ a\in A \, | \ \d^n (a)=0$
  for some natural $n\}$.

A $K$-derivation $\d $ of
 the algebra $A$ is a {\em locally nilpotent } derivation if for
 each element $a\in A$ there exists a natural number $n$ such
 that $\d^n(a)=0$. A $K$-derivation $\d $ is locally nilpotent iff
 $A=N(\d , A)$.

Given a ring $R$ and its derivation $d$. The {\em Ore extension}
$R[x;d]$ of $R$ is a ring freely generated over $R$ by $x$ subject
to the defining relations: $xr=rx+d(r)$ for all $r\in R$.
$R[x;d]=\oplus_{i\geq 0}Rx^i=\oplus_{i\geq 0}x^iR$ is a  left and
right free $R$-module. Given $r\in R$, a derivation $(\ad \,
r)(s):=[r,s]=rs-sr$ of $R$ is called an {\em inner} derivation of
$R$.

\begin{lemma}\label{dx=1}
 Let $A$ be an algebra over a field $K$ of characteristic zero and
$\d $ be a $K$-derivation of $A$ such that $\d (x)=1$ for some
$x\in A$. Then $N(\d ,A)=A^\d [x; d]$ is the Ore extension with
coefficients from the algebra $A^\d$, and the derivation $d$ of
the algebra $ A^\d$ is the restriction of the inner derivation
$\ad \, x $ of the algebra $A$ to its subalgebra $A^\d$. For each
$n\geq 0$, $N_n=\oplus_{i=0}^n\, A^\d x^i=\oplus_{i=0}^n\,
x^iA^\d$.
\end{lemma}

{\it Proof}. For each element $c\in C:=A^\d$,
$$ \d ([x,c])=[\d (x), c]+[x, \d (c)]=[1,c]+[x,0]=0,$$
thus $d (C)\subseteq C$, and $d $ is a $K$-derivation of the
algebra $C$.

First, we show that the $K$-subalgebra $N'$ of $N:=N(\d , A)$
generated by $C$ and $x$ is the Ore extension $C[x;d]$. We have
$N'=\sum_{i\geq 0}\, Cx^i$ since, for each $c\in C$, $xc-cx=d(c)$.
So, it remains to prove that the sum $\sum_{i\geq 0}\, Cx^i$ of
left $C$-modules is a direct sum. Suppose this is not the case,
then there is a nontrivial relation of degree $n>0$,
$$ c_0+c_1x+\cdots +c_nx^n=0, \;\; c_i\in C,\;\; c_n\neq 0.$$
We may assume that the degree $n$ of the relation above is the
{\em least} one. Then applying $\d $ to the relation above we
obtain the relation
$$ c_1+2c_2x+\cdots +nc_nx^{n-1}=0$$
of smaller degree $n-1$ since $nc_n\neq 0$ (char $K=0$), a
contradiction. So, $N'=C[x;d]$.

It remains to prove  that $N=N'$. The inclusion $N'\subseteq N$ is
obvious. In order to prove the inverse inclusion it suffices to
show that all subspaces $N_i$ belong to $N'$. We use induction on
$i$. The base of the induction is trivial since $N_0=C$. Suppose
that $i>0$, and $N_{i-1}\subseteq N'$. Let $u$ be an arbitrary
element of $N_i$. Then $\d (u)\in N_{i-1}\subseteq N'$. For an
arbitrary element $a=\sum \, c_jx^j\in N'$, we have $\d (b)=a$
where $b=\sum\, (j+1)^{-1}c_jx^{j+1}\in N'$. Therefore, in the
case of $a=\d (u)\in N'$, we have $\d (u)=\d (b)$ for some $b\in
N'$. Hence, $\d (u-b)=0$, and $u\in b+C\subseteq N'$. This means
that $N=N'$, as required. $\Box $

\begin{theorem}\label{8Nov05}
Let $A$ be an algebra over a field $K$ of characteristic zero, $\d
$ be a locally nilpotent $K$-derivation of the algebra $A$ such
that $\d (x)=1$ for some $x\in A$.  Then the $K$-linear  map $\phi
:=\sum_{i\geq 0} (-1)^i\frac{x^i}{i!}\d^i :A\ra A$ satisfies the
following properties:
\begin{enumerate}
\item  $\phi$ is a homomorphism of right $A^\d$-modules. \item
$\phi$ is a projection onto the algebra  $A^\d$:
$$ \phi : A=A^\d \oplus xA\ra A^\d \oplus xA, \;\; a+xb\mapsto a,
\;\; {\rm where}\;\; a\in A^\d, \; b\in A.$$ In particular,  ${\rm
im} (\phi )=A^\d$ and $\phi (y)=y$ for all $y\in A^\d$. \item
$\phi (x^i)=0$, $i\geq 1$. \item $\phi$ is an algebra homomorphism
provided $x\in Z(A)$, the centre of the algebra $A$.
\end{enumerate}
\end{theorem}

{\it Proof}. The map $\phi$ is well-defined since $\d$ is a
locally nilpotent derivation. It is obvious that $\phi$ is a
homomorphism of right $A^\d$-modules, $\phi (x)=0$, and $\phi
(y)=y$ for all $y\in A^\d$. An easy computation shows that $\d
\phi (z)=0$ for all $z\in A$, hence ${\rm im}(\phi )=A^\d$.

If $x$ is a {\em central} element then for $a,b\in A$:
\begin{eqnarray*}
\phi (ab)&=& \sum_{i\geq
0}(-1)^i\frac{x^i}{i!}\d^i(ab)=\sum_{i\geq
0}(-1)^i\frac{x^i}{i!}\sum_{j+k=i}{i\choose j}\d^j(a)\d^k(b)\\
&=& (\sum_{j\geq 0}(-1)^j\frac{x^j}{j!}\d^j(a))\, (\sum_{k\geq
0}(-1)^k\frac{x^k}{k!}\d^k(a))=\phi (a)\phi (b). \;\;
\end{eqnarray*}
For a {\em not} necessarily central $x$, repeating the above
computations we have, $\phi (x^i)=\phi (xx^{i-1})= \phi (x)\phi
(x^{i-1})=0$. Note that $A=\oplus_{i\geq 0}x^i A^\d $ (Lemma
\ref{dx=1}). Then $\phi$ is a projection onto $A^\d$,  $\phi :
a_0+xa_1+\cdots \mapsto a_0$, where $a_i\in A^\d$. $\Box$

The {\em Weyl } algebra $A_n=A_n(K)$ is a $K$-algebra generated by
$2n$ generators $q_1, \ldots , q_n$, $p_1, \ldots , p_n$ subject
to the defining relations:
$$ [p_i, q_j]=\d_{ij}, \;\; [p_i, p_j]=[q_i, q_j]=0\;\; {\rm
for\;\; all}\;\; i,j,$$ where $\d_{ij}$ is the Kronecker delta,
$[a,b]:=ab-ba$.

For the algebra $A_n\t P_m$ where $P_m=K[x_{2n+1}, \ldots
x_{2n+m}]$ is a polynomial algebra we use the following notation
$$ x_1:= q_1, \ldots , x_n:= q_n, x_{n+1}:= p_1, \ldots ,
x_{2n}:=p_n.$$ Then $A_n\t P_m=\bigoplus_{\alpha \in \mathbb{N}^s}
Kx^\alpha $ where $s:=2n+m$, $x^{\alpha }:= x_1^{\alpha_1 }\cdots
x_n^{\alpha_n }$, the {\em order} of the $x$'s is {\em fixed}. The
algebra $A_n\t P_m$ admits the finite set of {\em commuting
locally nilpotent} derivations, namely, the `partial derivatives':
$$ \der_1:= \frac{\der }{\der x_1}, \ldots , \der_s:= \frac{\der }{\der
x_s}.$$ Clearly, $\der_i=\ad (x_{n+i})$ and $\der_{n+i}=-\ad
(x_i)$, $i=1,\ldots , n$ (where $\ad (a):A\ra A$, $b\mapsto [a,b]$
is the {\em inner} derivation of the algebra $A$, $a\in A$).

For each $i=1, \ldots , s$, consider the maps from Theorem
\ref{8Nov05},
$$ \phi_i:=\sum_{k\geq 0}(-1)^k\frac{x_i^k}{k!}\der_i^k: A_n\t
P_m\ra A_n\t P_m.$$ For each $i=2n+1, \ldots ,s$, the map $\phi_i$
{\em commutes} with {\em all} the maps $\phi_j$. For each $i=1,
\ldots ,n$, the map $\phi_i$ commutes with all the maps $\phi_j$
but $\phi_{n+i}$, and the map $\phi_{n+i}$ commutes with all the
maps $\phi_j$ but $\phi_i$. Note that $A_n\t P_m=K\oplus V$ where
$V:=\bigoplus_{0\neq \alpha \in \mathbb{N}^s} Kx^\alpha$. Using
Theorem \ref{8Nov05}, we see that the map (the order is important)
\begin{equation}\label{phis}
\phi:=\phi_s\phi_{s-1}\cdots \phi_1:A_n\t P_m\ra A_n\t P_m, \;\;
a=\sum_{\alpha \in \mathbb{N}^s} \l_\alpha x^\alpha\mapsto \phi
(a)=\l_0,
\end{equation}
is a {\em projection} onto $K$.

The next result is a kind of a {\em Taylor} formula (note though
that even for polynomials this is {\em not} the Taylor formula.
The both formulae are essentially a formula for the {\em identity
map} and they give a presentation of an element as a series but
the formula below has one obvious advantage - it is `more
economical', i.e. there is no evaluation at $x=0$ as in the Taylor
formula).

\begin{theorem}\label{t8Nov05}
For any $a\in A_n\t P_m$,
$$ a=\sum_{\alpha \in \mathbb{N}^s}\phi (\frac{\der^\alpha}{\alpha
!} a)x^\alpha$$ where $s=2n+m$ and $\alpha !:=\alpha_1!\cdots
\alpha_s!$.
\end{theorem}

{\it Proof}. If $a=\sum \l_\alpha x^\alpha $, $\l_\alpha \in K$,
then, by (\ref{phis}), $\phi (\frac{\der^\alpha}{\alpha !}
a)=\l_\alpha$. $\Box $

So, the identity map ${\rm id} : A_n\t P_m \ra A_n\t P_m$ has a
nice presentation 
\begin{equation}\label{ida}
{\rm id}(\cdot ) = \sum_{\alpha \in \mathbb{N}^s}\phi
(\frac{\der^\alpha}{\alpha !} (\cdot ))x^\alpha .
\end{equation}

Let  $\Aut_K(A_n\t P_m)$ be the group of $K$-algebra automorphisms
of the algebra $A_n\t P_m$. Given an automorphism $\s
\in\Aut_K(A_n\t P_m)$. It is uniquely determined by the elements
\begin{equation}\label{xsh1}
x_1':= \s (x_1), \ldots , x_s':=\s (x_s)
\end{equation}
of the algebra $A_n\t P_m$. The centre $Z:=Z(A_n\t P_m)$ of the
algebra $A_n\t P_m$ is equal to $P_m$. Therefore, the restriction
$\s|_{P_m}\in \Aut_K(P_m)$, and so
$$ \D :=\det (\frac{\der x_{2n+i}'}{\der x_{2n+j}})\in K^*$$
where $i,j=1, \ldots , n$. The corresponding (to the elements
$x_1',\ldots , x_s'$) `partial derivatives' (the set of commuting
locally nilpotent derivations of the algebra $A_n\t P_m$)
\begin{equation}\label{xsh2}
\der_1':= \frac{\der }{\der x_1'}, \ldots , \der_s':= \frac{\der
}{\der x_s'}
\end{equation}
are equal to 
\begin{equation}\label{dad1}
\der_i':= \ad (\s (x_{n+i})), \;\; \der_{n+i}':= -\ad (\s
(x_{i})), \;\; i=1, \ldots , n,
\end{equation}

\begin{equation}\label{dad2}
\der_{2n+j}'  := \D^{-1} \det
 \begin{pmatrix}
  \frac{\der \s (x_{2n+1})}{\der x_{2n+1}} & \cdots & \frac{\der \s (x_{2n+1})}{\der x_{2n+m}} \\
  \vdots & \vdots & \vdots \\
\frac{\der }{\der x_{2n+1}} & \cdots & \frac{\der }{\der x_{2n+m}}\\
 \vdots & \vdots & \vdots \\
\frac{\der \s (x_{2n+m})}{\der x_{2n+1}} & \cdots & \frac{\der \s (x_{2n+m})}{\der x_{2n+m}} \\
\end{pmatrix}, \;\;\; j=1, \ldots , m,
\end{equation}
where we `drop' $\s (x_{2n+j})$ in the determinant $\det
(\frac{\der \s (x_{2n+k})}{\der x_{2n+l}})$.

For each $i=1, \ldots , s$, let 
\begin{equation}\label{dad3}
\phi_i':= \sum_{k\geq 0}(-1)^k\frac{(x_i')^k}{k!}(\der_i')^k:
A_n\t P_m\ra A_n\t P_m
\end{equation}
and (the order is important) 
\begin{equation}\label{dad4}
\phi_\s := \phi_s'\phi_{s-1}'\cdots \phi_1'.
\end{equation}

\begin{theorem}\label{i8Nov05}
{\rm (The Inversion Formula)} For each $\s \in \Aut_K(A_n\t P_m)$
 and $a\in A_n\t P_m$,
 $$ \s^{-1}(a)=\sum_{\alpha \in \mathbb{N}^s}\phi_\s
 (\frac{(\der')^\alpha}{\alpha!}a)x^\alpha , $$
 where $(\der')^{\alpha} :=(\der_1')^{\alpha_1}\cdots
 (\der_s')^{\alpha_s}$ and  $s=2n+m$.
\end{theorem}

{\it Proof}. By Theorem \ref{t8Nov05}, $a=\sum_{\alpha \in
\mathbb{N}^s}\phi_\s
 (\frac{(\der')^\alpha}{\alpha!}a)(x')^\alpha$. Applying $\s^{-1}$
 we have the result
 $$ \s^{-1}(a)=\sum_{\alpha \in
\mathbb{N}^s}\phi_\s
 (\frac{(\der')^\alpha}{\alpha!}a)\s^{-1}((x')^\alpha ) =\sum_{\alpha \in
\mathbb{N}^s}\phi_\s
 (\frac{(\der')^\alpha}{\alpha!}a)x^\alpha . \;\;\; \Box $$

\begin{corollary}\label{ci8Nov05}
The question in the Abstract has an affirmative answer iff both
the Jacobian conjecture and the Dixmier problem have (in more
detail, $JD_n$ $\Leftrightarrow $ $JC_n+ DP_n$).
\end{corollary}

{\it Proof}. $ (\Rightarrow )$ Obvious (the JC and the DP are
special cases).

$ (\Leftarrow)$ Suppose that a $K$-algebra endomorphism $\s :A_n\t
P_m\ra A_n\t P_m$ satisfies the condition $\det (\frac{\der \s
(x_i)}{\der x_j})=1$. Note that $P_m$ is the {\em centre} of the
algebra $A_n\t P_m$. Then the JC implies $\s |_{P_m}\in
\Aut_K(P_m)$. Without loss of generality one can assume that $\s
|_{P_m}={\rm id}$. Let $Q_m$ be the field of fractions of $P_m$.
Then $\s$ can be extended to an endomorphism of the algebra $A_n\t
Q_m$. By the DP, $\s \in \Aut_K(A_n\t Q_m)$. By Theorem
\ref{i8Nov05}, $\s^{-1}(A_n\t P_m)\subseteq A_n\t P_m$, and so $\s
\in \Aut_K(A_n\t P_m)$. $\Box $

{\it Remark}. Note that an algebra endomorphism $\s $ of the
algebra $A_n\t P_m$ satisfying $\det (\frac{\der \s (x_i)}{\der
x_j})\in K^*$ is automatically an algebra {\em monomorphism}: $\s
|_{P_m}$ is an algebra monomorphism, it induces an algebra
monomorphism, say $\s $,  on the field of fractions $Q_m$ of
$P_m$, hence $\s$ can be extended to an algebra endomorphism of
the {\em simple} algebra $A_n\t Q_m$, hence $\s$ is an algebra
{\em monomorphism}.


\section{The Degree of Inverse
Automorphism}\label{DEGINVA}

In this section, Theorem \ref{9Nov05} and Corollary \ref{c9Nov05}
are proved.

For an automorphisms $\s$ and $\s^{-1}$ we keep the notation from
the previous sections.

An automorphism $\s \in \Aut_K(A_n\t P_m)$ is uniquely determined
by  $\s (x_1), \ldots , \s (x_s)$. The {\em degree} of the
automorphism $\s $ is defined as
$$ \deg \, \s :=\max \{ \deg \, \s (x_i)\, | \, i=1, \ldots ,
s\}$$
 where the {\em degree} $\deg \, a$ of an element $a=\sum_{\alpha
 \in \mathbb{N}^s} \l_\alpha x^\alpha\in A_n\t P_m$ is defined as
 $$ \deg \, a := \max \{ |\alpha | := \alpha_1+\cdots +\alpha_s\,
 | \, \l_\alpha \neq 0\}.$$

\begin{theorem}\label{BCRRM}
\cite{BCW}, \cite{Rus-Win84}. Given $\s \in \Aut_K(P_m)$. Then
$\deg \, \s^{-1}\leq (\deg \, \s )^{m-1}$.
\end{theorem}

By Theorem \ref{i8Nov05}, $\s^{-1}(x_i)=\sum_{\alpha \in
\mathbb{N}^s}\phi_\s (\frac{(\der')^\alpha}{\alpha !}x_i)x^\alpha
$, then applying $\s $ to this equality and using the fact that
$\s (x^\alpha )=(x')^\alpha$, we have the equality

\begin{equation}\label{dad5}
\s^{-1}(x_i')=\s^{-1}\s (x_i)=\s \s^{-1} (x_i)=\sum_{\alpha \in
\mathbb{N}^s}\phi_\s (\frac{(\der')^\alpha}{\alpha
!}x_i)(x')^\alpha .
\end{equation}
The next lemma follows directly from (\ref{dad5}) and it gives the
{\em exact} value for the degree of $\s^{-1}$.

\begin{lemma}\label{n9Nov05}
Given $\s \in \Aut_K(A_n\t P_m)$. Then
$$ \deg \, \s^{-1} = \max \{ \deg' (x_i)\, | \, i=1, \ldots ,
s\}$$
 where, for $a=\sum_{\alpha \in \mathbb{N}^s}
 \l_\alpha'(x')^\alpha \in A_n\t P_m$, $\deg' (a):= \max \{
 |\alpha | \, | \, \l_\alpha'\neq 0\}$.
\end{lemma}

\begin{theorem}\label{9Nov05}
Given $\s \in \Aut_K(A_n\t P_m)$. Then
$$ \deg \, \s^{-1}\leq (\deg \, \s )^{s-1}$$
where $s:=2n+m$.
\end{theorem}

{\it Proof}.  The algebra $A:=A_n\t P_m=\cup_{\alpha \in
\mathbb{N}^s} Kx^\alpha =\cup_{i\geq 0}\CA_i$ is a {\em filtered}
algebra ($\CA_i \CA_j\subseteq \CA_{i+j}$ for all $i,j\geq 0$)
where $\CA_i := \oplus_{|\alpha |\leq i}Kx^\alpha $ and
$$ \dim_K(\CA_i)={i+s\choose s}=\frac{(i+s)(i+s-1)\cdots (i+1)}{s!}=\frac{i^s}{s!}+\cdots , \;\; i\geq 0,$$
where here and everywhere the three dots mean `smaller' terms. The
filtration $\{ \CA_i\}$ is a {\em standard} filtration (we use the
terminology of the book of H. Krause and T. Lenagan, \cite{KL},
where one can find all the missed definitions), so the {\em
Gelfand-Kirillov} dimension  of the algebra $A$ is $\GK (A)=s$.
The associative graded algebra ${\rm gr} \, A:= \oplus_{i\geq 0}
\, \CA_i/\CA_{i-1}$ is canonically isomorphic to a polynomial
algebra in $s$ variables. So, the algebra $A$ is an {\em almost
commutative} algebra. Given a finitely generated $A$-module
$M=AM_0=\cup_{i\geq 0}M_i$, $M_i:=\CA_i M_0$, where $M_0$ is a
finite dimensional generating space for the module $M$, then there
exists a polynomial (so-called, the {\em Hilbert polynomial} of
$M$) $p_M\in \mathbb{Q}[t]$ such that
$$ \dim_K(M_i)=p_M(i)=\frac{e(M)i^{\GK (M)}}{\GK (M)!}+\cdots ,
\;\; i\gg  0,$$ where $e(M)\in \mathbb{N}$ is the {\em
multiplicity} of $M$. {\em All} algebras involved in this proof
will be {\em algebras of the type} $A_k\t P_l$, so we will use
Hilbert polynomials and multiplicity  for certain modules.

Fix $\nu \in \{ 1, \ldots , s\}$. Then 
\begin{equation}\label{Ldn1}
A=\oplus_{k\geq 0}\L (x_\nu')^k, \;\; {\rm where}\;\;\; \L =\L
(\nu ):= K\langle x_1', \ldots , \widehat{x_\nu '}, \ldots ,
x_s'\rangle
\end{equation}
is an algebra of type $A_k\t P_l$ (the hat over the symbol means
that it is missed), and $\GK (\L )=s-1$. A nonzero element $a\in
A$ is a {\em unique} sum
$$ a=a_0+a_1x_\nu'+\cdots +a_d(x_\nu')^d, \;\; a_i\in \L, \;\;
a_d\neq 0.$$ Then $d$ is called the $x_\nu'$-{\em degree} of the
element $a$ denoted $\deg_{x_\nu'} (a)$.

Fix $ j\in \{ 1, \ldots , s\}$. Let $d_j:= \deg_{x_\nu'} (x_j)$.
By Lemma \ref{n9Nov05}, in order to finish the proof of this
theorem it suffices to show that $d_j\leq (\deg \, \s )^{s-1}$
(the field $K$ has characteristic zero, in particular it is {\em
infinite}, so making suitable `linear changes of variables' $\{
x_\mu'\}$ (i.e. up to a  {\em linear algebra automorphism} of
$A_n\t P_m$) one can assume that $d_j=\deg \, \s^{-1}$).

Clearly,
$$ A\supseteq M\oplus Ax_j, \;\; M:=\oplus_{k=0}^{d_j-1} \L
(x_\nu')^k.$$ Then $M$ can be seen as a $\L$-{\em submodule} of
$A/Ax_j$. Consider the filtration $\{ \CP_i\}$  on the $A$-module
$A/Ax_j$ induced by the $1$-dimensional generating space
$K\overline{1}$, $\overline{1}:=1+Ax_j$. Then
$$ \dim_K(\CP_i)={i+s-1\choose s-1} =
\frac{i^{s-1}}{(s-1)!}+\cdots , \;\; i\geq 0. $$ Let $\{ \CP_i'\}$
be the {\em standard} filtration of the algebra $\L$ (with respect
to the generating set $x_1', \ldots , \widehat{x_\nu'}, \ldots ,
x_s'$). Clearly,

$$ \CP_i'\subseteq \CP_{i(\deg \, \s)}, \;\; i\geq 0.$$
Fix a natural number, say $t$, such that $(x_\nu')^k\in \CP_t$ for
all $k=0, \ldots , d_j-1$. Then, for all $i\gg 0$, $M_i:=
\oplus_{k=0}^{d_j-1}\CP_i'(x_\nu')^k\subseteq \CP_{i(\deg \,
\s)+t}$. Therefore,
\begin{eqnarray*}
\dim_K(\CP_{i(\deg\, \s )+t}) & = & \frac{(\deg \, \s )^{s-1} i^{s-1} }{(s-1)!}+\cdots \geq \dim_K(M_i) \\
 &\geq & \sum_{k=0}^{d_j-1}{i+s-1-t\choose s-1}=
 \frac{d_ji^{s-1} }{ (s-1)! }+\cdots , \;\; i\gg 0.
\end{eqnarray*}
Hence, $d_j\leq (\deg \, \s )^{s-1}$, as required.  $\Box $

Recall that a {\em commutative} ring $R$ is called {\em reduced}
iff its
 {\em nil-radical} is equal to zero ($\gn (R):=\cap_{\gp \in \Spec
 (R)} \gp =0$) iff  a zero is the only nilpotent element of $R$.

\begin{corollary}\label{c9Nov05}
If a commutative ring (not necessarily a field)  $K$ is a reduced
$\mathbb{Q}$-algebra and $\s \in \Aut_K(A_n\t P_m)$. Then $$\deg\,
\s^{-1}\leq (\deg \, \s )^{s-1}.$$
\end{corollary}

{\it Proof}. We write $A_n(K)$ and $P_m(K)$ to indicate the base
ring $K$. If $K$ is a {\em domain}  with a field of fractions $F$
then $\s $ can be extended to an element of
$\Aut_F(A_n(F)\t_FP_m(F))$ and the result follows from Theorem
\ref{9Nov05}.

In the general situation, by the previous case, for each prime
ideal $\gp $ of $K$,  reduction modulo $\gp$ gives an element
$\s_\gp \in \Aut_{K/\gp }(A_n(K/\gp )\t_{K/\gp}P_m(K/\gp ))$.
Since $K/\gp$ is a domain,
$$ \deg \, \s_\gp^{-1}\leq (\deg \, \s_\gp )^{s-1}\leq (\deg \, \s
)^{s-1} \;\; {\rm for \;\; all}\;\; \gp \in \Spec (K),$$ which
implies $\deg\, \s^{-1}\leq (\deg \, \s )^{s-1}$ since $K$ is
reduced. $\Box $

{\it Remark}. In the polynomial case, $P_m$, if $K$ is {\em not
reduced} then there is no uniform upper bound for the degree $\deg
\, \s^{-1}$ depending only on $m$ and $\deg \, \s$, see p.56--57,
\cite{vdEssenPA}: let $P_1=K[x]$ and $K:=\mathbb{Q}[T]/(T^l)$
where $\s :x\mapsto x':= x-tx^2$ where $t:=T+(T^l)\in K$. Then
$\deg_x(x')\geq \frac{l-1}{2}+1$ (p.57, \cite{vdEssenPA}). {\em
The same is true for automorphisms of the algebras} $A_n\t P_m$
since the automorphism $\s$ can extended to a $K$-automorphism of
the 1st Weyl algebra $A_1(K)$ by the rule
$$ \frac{d}{dx}\mapsto \frac{d}{dx'}=\frac{dx}{dx'}\frac{d}{dx}=
(\frac{dx'}{dx})^{-1}\frac{d}{dx}=\frac{1}{1-2tx}\frac{d}{dx}=
(1+2tx+\cdots +(2tx)^{m-1})\frac{d}{dx}.$$


\section{The inversion formula for automorphism of a power series
algebra}\label{IFFPSA}

In this section,  the  inversion formula for an automorphism of
the power series algebra is obtained (Theorem \ref{hPi8Nov05}).

\begin{lemma}\label{jhP8Nov05}
Let $A$ be an algebra over a field $K$ of characteristic zero, $\d
$ be a  $K$-derivation of the algebra $A$ (not necessarily locally
nilpotent) such that $\d (x)=1$ for a central element  $x\in A$.
Suppose that the algebra $A$ is complete in $\gm $-adic topology
and $x\in \gm$ ($\gm $ is a right ideal of $A$). Then the
$K$-linear map $\phi :=\sum_{i\geq 0} (-1)^i\frac{x^i}{i!}\d^i
:A\ra A$ satisfies the following properties:
\begin{enumerate}
\item $\phi (x)=0$.
 \item $\phi$ is an algebra homomorphism of $A$.
 \item  $\phi$ is a homomorphism of left/right
$A^\d$-modules. \item 
${\rm im} (\phi )=A^\d$ and $\phi (y)=y$ for all $y\in A^\d$.
\end{enumerate}
\end{lemma}

{\it Proof}. The map $\phi$ is well-defined. Then the proof is a
repetition of the proof of Theorem \ref{8Nov05}. $\Box$

{\it Remark}. If for an arbitrary $K$-algebra $A$ the infinite sum
$\phi :=\sum_{i\geq 0} (-1)^i\frac{x^i}{i!}\d^i$  makes sense then
Lemma \ref{jhP8Nov05} holds.

Let $\widehat{P}_n:= K[[x_1, \ldots , x_m]]$ be {\em an algebra of
formal power series} in $x_i$ (= the completion of the polynomial
algebra $P_m$ at the maximal ideal $\gm =(x_1, \ldots , x_m)$).
The partial derivatives $\der_1, \ldots , \der_m$ are a set of
commuting continuous (in $\gm$-adic topology) derivations of the
algebra $\widehat{P}_m$ satisfying $\der_i(x_j)=\d_{ij}$.

Using Lemma \ref{jhP8Nov05}, we have a set of {\em commuting}
algebra endomorphisms:
$$ \phi_i:=\sum_{k\geq 0}(-1)^k\frac{x_i^k}{k!}\der_i^k:
\widehat{P}_m\ra\widehat{P}_m, \;\; i=1, \ldots , m.$$

Each $\phi_i$ is a {\em projection} onto ${\rm im} \,
\phi_i=K[[x_1, \ldots , \widehat{x_i}, \ldots , x_n]]$,
$\phi_i(x_i)=0$, and $\phi_i:\widehat{P}_m\ra \widehat{P}_m$ is a
homomorphism of $K[[x_1, \ldots , \widehat{x_i}, \ldots ,
x_n]]$-modules ($\widehat{x}_i$ means that $x_i$ is missed).

The algebra endomorphism  
\begin{equation}\label{hPphis}
\phi:=\phi_1\cdots \phi_m:\widehat{P}_m=K\oplus \widehat{P}_m\gm
\ra \widehat{P}_m=K\oplus \widehat{P}_m\gm, \;\; a=\sum_{\alpha
\in \mathbb{N}^m} \l_\alpha x^\alpha\mapsto \phi (a)=\l_0,
\end{equation}
is a {\em projection} onto $K$.

\begin{theorem}\label{hPt8Nov05}
For any $a\in \widehat{P}_m$,
$$ a=\sum_{\alpha \in \mathbb{N}^m}\phi (\frac{\der^\alpha}{\alpha
!} a)x^\alpha .$$
\end{theorem}

{\it Proof}. If $a=\sum \l_\alpha x^\alpha $, $\l_\alpha \in K$,
then, by (\ref{hPphis}), $\phi (\frac{\der^\alpha}{\alpha !}
a)=\l_\alpha$. $\Box $

So, the identity map ${\rm id} : \widehat{P}_m \ra \widehat{P}_m$
has a nice presentation 
\begin{equation}\label{ida}
{\rm id}(\cdot ) = \sum_{\alpha \in \mathbb{N}^m}\phi
(\frac{\der^\alpha}{\alpha !} (\cdot ))x^\alpha .
\end{equation}

Let  $\Aut_K(\widehat{P}_m)$ be the group of continuous
$K$-algebra automorphisms of the algebra $\widehat{P}_m$. Given an
automorphism $\s \in\Aut_K(\widehat{P}_m)$. It is uniquely
determined by the elements 
\begin{equation}\label{hPxsh1}
x_1':= \s (x_1), \ldots , x_m':=\s (x_m)
\end{equation}
of the algebra $\widehat{P}_m$ (all the
$x_i'\in\widehat{P}_m\gm$). Then $ \D :=\det (\frac{\der
x_i'}{\der x_j})\in (\widehat{P}_m)^*$, the group of all
invertible elements of $\widehat{P}_m$. The corresponding (to the
elements $x_1',\ldots , x_m'$) `partial derivatives' (the set of
commuting continuous derivations of the algebra $\widehat{P}_m$)
\begin{equation}\label{hPxsh2}
\der_1':= \frac{\der }{\der x_1'}, \ldots , \der_m':= \frac{\der
}{\der x_m'}
\end{equation}
are equal to 
\begin{equation}\label{hPdad2}
\der_i'  := \D^{-1} \det
 \begin{pmatrix}
  \frac{\der \s (x_1)}{\der x_1} & \cdots & \frac{\der \s (x_1)}{\der x_m} \\
  \vdots & \vdots & \vdots \\
\frac{\der }{\der x_1} & \cdots & \frac{\der }{\der x_m}\\
 \vdots & \vdots & \vdots \\
\frac{\der \s (x_m)}{\der x_1} & \cdots & \frac{\der \s (x_m)}{\der x_m} \\
\end{pmatrix}, \;\;\; j=1, \ldots , m,
\end{equation}
where we `drop' $\s (x_i)$ in the determinant $\det (\frac{\der \s
(x_i)}{\der x_j})$.

For each $i=1, \ldots , m$, let 
\begin{equation}\label{hPdad3}
\phi_i':= \sum_{k\geq 0}(-1)^k\frac{(x_i')^k}{k!}(\der_i')^k:
\widehat{P}_m \ra \widehat{P}_m
\end{equation}
and 
\begin{equation}\label{hPdad4}
\phi_\s := \phi_1'\cdots \phi_m'.
\end{equation}

\begin{theorem}\label{hPi8Nov05}
{\rm (The Inversion Formula)} For each $\s \in
\Aut_K(\widehat{P}_m)$
 and $a\in \widehat{P}_m$,
 $$ \s^{-1}(a)=\sum_{\alpha \in \mathbb{N}^m}\phi_\s
 (\frac{(\der')^\alpha}{\alpha!}a)x^\alpha . $$
\end{theorem}

{\it Proof}. By Theorem \ref{hPt8Nov05}, $a=\sum_{\alpha \in
\mathbb{N}^m}\phi_\s
 (\frac{(\der')^\alpha}{\alpha!}a)(x')^\alpha$. Applying $\s^{-1}$
 we have the result
 $$ \s^{-1}(a)=\sum_{\alpha \in
\mathbb{N}^m}\phi_\s
 (\frac{(\der')^\alpha}{\alpha!}a)\s^{-1}((x')^\alpha ) =\sum_{\alpha \in
\mathbb{N}^m}\phi_\s
 (\frac{(\der')^\alpha}{\alpha!}a)x^\alpha . \;\;\; \Box $$

\begin{corollary}\label{1hP8Nov05}
Let $\s$ be an algebra endomorphism  of the polynomial algebra
$P_m=K[x_1, \ldots , x_m]$ satisfying $\det (\frac{\der \s
(x_i)}{\der x_j})\in K^*$ and $\s (\gm )\subseteq \gm $ where $\gm
:=(x_1, \ldots , x_n)$. Then
\begin{enumerate}
\item the algebra endomorphism $\phi_\s :P_m=K\oplus \gm \ra
P_m=K\oplus \gm$, $\l +\sum x_ia_i\mapsto \l$, $\l \in K$, $a_i\in
P_m$ (see (\ref{hPdad4})) is a projection onto $K$. \item
$\cap_{i=1}^m\ker_{P_m}(\der_i')=\cap_{i=1}^m\ker_{\widehat{P}_m}(\der_i')=K$.
\item $\cap_{i=1}^m N(\der_i', P_m)=\cap_{i=1}^m N(\der_i',
\widehat{P}_m)=\s (P_m)$.
\end{enumerate}

\end{corollary}

{\it Proof}. $1$. The two conditions guarantee that the extension
of the $\s$ to a continuous (in $\gm $-adic topology) algebra
endomorphism, say $\s$,  of $\widehat{P}_m$ is, in fact, an
automorphism. By (\ref{hPphis}), the endomorphism $\phi_\s :
\widehat{P}_m\ra \widehat{P}_m$ is the projection onto $K$, hence
its restriction $\phi_\s : P_m\ra P_m$ is a projection onto $K$ as
well.

$2$. $K\subseteq \cap_{i=1}^m\ker_{P_m}(\der_i') \subseteq
\cap_{i=1}^m\ker_{\widehat{P}_m}(\der_i') =\phi_\s
(\widehat{P}_m)=K$.

$3$. Note that $\der_1',\ldots , \der_m'$ is a set of commuting
{\em locally nilpotent} derivations of the algebra $\cap_{i=1}^m
N(\der_i', P_m)$. By Corollary \ref{1genbassdx=1}, $\cap_{i=1}^m
N(\der_i', P_m)=C[\s (x_1), \ldots, \s (x_m)]$ and, similarly,
$\cap_{i=1}^m N(\der_i', \widehat{P}_m)=C[\s (x_1), \ldots, \s
(x_m)]$ where $C=\cap_{i=1}^m \ker_{P_m}(\der_i')=\cap_{i=1}^m
\ker_{\widehat{P}_m}(\der_i')=K$, by statement $2$.  $\Box $


\section{A Characterization of the Weyl Algebras and the
Polynomial Algebras in terms of commuting locally nilpotent
derivations}

 The Weyl algebras and the polynomial algebras are, in some sense,
 special algebras in the class of all algebras. Theorem \ref{Wpolchar}
 explains this fact in
 terms of commuting locally nilpotent derivations.

We say that derivations $\d_1, \ldots , \d_s$ of an algebra $A$
have {\em generic kernels} iff the intersections $\cap_{i\in
I}\ker \, \d_i$, $I\subseteq \{ 1, \ldots , s\}$,  are {\em
distinct} (iff $\cap_{j\neq i}\ker \, \d_j$, $i=1, \ldots , s$ are
distinct).

\begin{lemma}\label{dix=1}
Let $A$ be an algebra over a field $K$ of characteristic zero, and
$\d_1, \ldots , \d_s$ be commuting locally nilpotent
$K$-derivations of $A$ that have generic kernels and with $\Gamma
:=\cap_{i=1}^s\, \ker \, \d_i$, a division ring such that $\Gamma
\neq A$. Then there exist nonzero elements $x_i\in
C_i':=\cap_{j\neq i}\ker\, \d_j$, $i=1, \ldots , s$ such that
$\d_i(x_j)=\d_{ij}$, the Kronecker delta.
\end{lemma}

{\it Proof}. Let us consider first the case when  $s=1$. The
derivation $\d_1$ is locally nilpotent and nonzero since $\ker \,
\d =\Gamma \neq A$. So we can find an element $y\in A$ such that
$0\neq \l :=\d_1(y)\in \ker \, \d_1=\Gamma$. Then $\d_1(x)=1$
 for $x=\l^{-1}y$.

Suppose now that $s>1$. For each $i=1, \ldots , s$, let
$\overline{\d_i}$ be the restriction of the derivation $\d_i$ to
the subalgebra $C_i'$ of $A$. The kernel of $\overline{\d_i}$ is
equal to $\Gamma$, so by the previous argument one can find  an
element, say $x_i\in C_i'$, satisfying $\d_i(x_i)=1$. Obviously,
$\d_i(x_j)=\d_{ij}$. $\Box $

Let us recall one of the key results of symplectic algebra.

\begin{lemma}\label{abfcanb}
Let $\Phi $ be an antisymmetric bilinear form on a finite
dimensional vector space $V$, and let $y_1, \ldots , y_m$ be a
basis of the kernel of $\Phi $. Then we can complete the set $y_1,
\ldots , y_m$ to the basis $p_1, \ldots , p_n, q_1, \ldots , q_n,
y_1, \ldots , y_l$ of $V$ such that $\Phi(p_i, q_j)=\d_{ij}$,
$\Phi (p_i, p_j)=\Phi (q_i, q_j)=0$ for all $i,j$.
\end{lemma}

The theorem below gives a characterization of the algebras of the
type $A_n\t P_m$ (the Weyl algebras with polynomial coefficients)
in terms of commuting locally nilpotent derivations.

\begin{theorem}\label{Wpolchar}
Let $A$ be an algebra over a field $K$ of characteristic zero.
Then the following statements  are equivalent.
\begin{enumerate}
\item There exist commuting locally nilpotent nonzero
$K$-derivations $\d_1, \ldots , \d_s$ of the algebra $A$ with
generic kernels $C_i=\ker \, \d_i$ satisfying $\cap_{i=1}^s\,
C_i=K$. \item The algebra $A$ is an iterated Ore extension
$$ A=K[x_1][x_2;d_2][x_3;d_3]\cdots [x_s;d_s]$$
such that $\l_{ij}:=d_i(x_j)\in K$ for all $i>j$, and
$\d_i(x_j)=\d_{ij}$, the Kronecker delta. \item The algebra $A$ is
isomorphic to the tensor product $A_n\t P_m$ (over $K$) of the
Weyl algebra $A_n$     with a polynomial algebra $P_m$ in $m$
indeterminates, and $2n+m=s$.


Suppose that the (equivalent) conditions above hold. Then
\begin{enumerate}
\item  the elements $x_1, \ldots , x_s$ are uniquely determined up
to scalar addition. \item  $n= \frac{1}{2} {\rm rk} (\L )$ and
$m=s-2n=\dim \, \ker (\L ) $ where $\L =(\l_{ij})$ is the
antisymmetric $s\times s$ matrix with lower diagonal entries
$\l_{ij}$ as above. \item  For each $i$, the algebra $C_i$ is an
iterated Ore extension
$$ K[x_1][x_2;d_2]\cdots [x_{i-1};d_{i-1}][x_{i+1};d_{i+1}]\cdots [x_s;d_s]$$
with $d_i(x_j)=\l_{ij}$ as above. Hence, $C_i\simeq A_{m_i}\t
P_{l_i}$ with $2m_i+l_i=s-1$. \item  The Gelfand-Kirillov
dimension $\GK (A)=s$.
\end{enumerate}\end{enumerate}
\end{theorem}

{\it Remark}. As an abstract algebra the iterated Ore extension
$A$ from the second statement is generated by the elements $x_1,
\ldots x_s$ subject to the defining relations
$$ x_ix_j-x_jx_i=\l_{ij}, \;\; {\rm for \;\; all}\;\; i>j.$$
So, for any permutation $i_1, \ldots , i_s$, of the indices $1,
\ldots , s$, the algebra $A$ is the iterated Ore extension
$$ K[x_{i_1}][x_{i_2}; d_{i_2}]\cdots [x_{i_s}; d_{i_s}]$$
with $d_\alpha (x_\beta )=\l_{\alpha \beta}$, if $\alpha >\beta $,
and $-\l_{\alpha \beta}$, if $\alpha <\beta $.

{\it Proof}. $(1\Rightarrow 2)$ We use induction on $s$. Let
$s=1$. Since $\d_1$ is a locally nilpotent nonzero derivation with
$\ker \, \d_1=K$ we can find an element $x\in A$ such that $\d_1
(x)=1$. By Lemma \ref{dx=1}, $A=K[x]$ is a polynomial algebra.

Suppose that $s>1$ and the result is true for $s-1$. By Lemma
\ref{dix=1}, we can find  a nonzero element $x_s\in A$ such that $
\d_i(x_s)=\d_{is}$.  By Lemma \ref{dx=1}, $A=N(\d_s,A)=C_s[x_s;
d_s]$ is an Ore extension with coefficients from $C_s$ where the
derivation $d_s$ of $C_s$  is the restriction of the inner
derivation $\ad \, x_s$ of $A$ to $C_s$. The derivations $\d_i$
commute, thus $\d_i (C_s)\subseteq C_s$ for all $i$.
 We denote by
$\overline{\d_1},\ldots , \overline{\d_{s-1}}$ the restrictions of
the derivations $\d_1, \ldots , \d_{s-1}$ to the subalgebra $C_s$
of $A$. Then $\overline{\d_1},\ldots , \overline{\d_{s-1}}$ are
commuting locally nilpotent nonzero $K$-derivations of the algebra
$C_s$ with generic kernels $\overline{C_1}, \ldots ,
\overline{C_{s-1}}$ satisfying $\cap_{i=1}^{s-1}\,
\overline{C_i}=\cap_{i=1}^s\, C_i=K$. By induction, the algebra
$\overline{C_s}$ is an iterated Ore extension
$$ K[x_1][x_2;d_2]\cdots  [x_{s-1};d_{s-1}]$$
 with $d_i(x_j)=\l_{ij}\in K$ for all $i,j$ less than $s$. It
 remains to prove that the elements
 $\l_{si}:=d_s(x_i)=[x_s,x_i]\in C_s$ belong to $K$, for all
 $i<s$. For $i,j<s$,
 $$
 \d_j(\l_{si})=\d_j([x_s,x_i])=[\d_j(x_s),x_i]+[x_s,\d_j(x_i)]=[0,x_i]+[x_s,
 \l_{ji}]=0,$$
 thus $\l_{si}\in \cap_{k=1}^s\, C_k=K$, as required. Let $x_1',
 \ldots , x_s'$ be elements of $A$ satisfying
 $\d_i(x_j')=\d_{ij}$. Then $\d_i (x_j-x_j')=0$, so $x_j-x_j'\in
 \cap_{k=1}^s \, C_k=K$, this proves $(a)$.

$(2\Rightarrow 3)$ The $s$-dimensional vector subspace
$V=Kx_1\oplus \cdots \oplus Kx_s$ of the iterated Ore extension
$A$ as in statement 2 is equipped with the antisymmetric bilinear
form:
$$ V\times V\ra K,\;\; (u,v)\ra [u,v]:=uv-vu.$$
Then $\L =(\l_{ij})$ is the matrix of this form in the basis $x_1,
\ldots , x_s$ of $V$. By Lemma \ref{abfcanb}, we can choose a
basis $p_1, \ldots , p_n, q_1, \ldots , q_n, y_1, \ldots , y_m$ of
the vector space $V$ such that
$$ [p_i, q_j]=\d_{ij}, \;\; [p_i, p_j]=[q_i,q_j]=0, $$
$$ [p_i,y_k]=[q_i, y_k]=[y_k, y_{k'}]=0,$$
for all possible $i$, $j$, $k$, and $k'$. So, $A=A_n\t P_m$ where
$ A_n=K[p_1, \ldots , p_n, q_1,\ldots , q_n]$ is the Weyl algebra
and $P_m=K[y_1, \ldots , y_m]$ is a polynomial algebra in $m$
indeterminates (see the Remark above).

Clearly, $2n+m=s$, $ n={1\over 2}{\rm rk} (\L )$, and $
m=s-2n=\dim \, \ker (\L ),$ this proves $(b)$.

$(3\Rightarrow 1)$ If $A=A_n\t P_m$ (as above) then
$$ \ad \, p_1, \ldots , \ad \, p_n,\,  \ad \, q_1, \ldots , \ad \,
q_n,\,  \frac{\der}{\der y_1}, \ldots , \frac{\der }{ \der y_m}$$
are commuting locally nilpotent nonzero $K$-derivations of the
algebra $A$ with generic kernels (recall that, for each element
$a\in A$, $\ad \, a:A\ra A$, $x\mapsto [a,x]:=ax-xa$, is the {\em
inner derivation} of the algebra $A$)
$$ \ker \, \ad \, p_i=K[p_i]\t K[p_1, \ldots , \widehat{p_i},
\ldots , p_n, q_1, \ldots , \widehat{q_i}, \ldots , q_n]\t
P_m\simeq A_{n-1}\t P_{m+1},$$
$$ \ker \, \ad \, q_i=K[q_i]\t K[p_1, \ldots , \widehat{p_i},
\ldots , p_n, q_1, \ldots , \widehat{q_i}, \ldots , q_n]\t
P_m\simeq A_{n-1}\t P_{m+1},$$
$$ \ker \, \frac{\der }{ \der y_j}=A_n\t K[y_1, \ldots , \widehat{y_j},
\ldots , y_m]\simeq A_n\t P_{m-1},$$ such that their intersection
is $K$ (where  hat over a symbol means that it is missed in the
list).  The statement $(c)$ is evident now.  Finally, the
Gelfand-Kirillov dimension of the algebra $A$ over the field $K$
is equal to
$$ \GK (A)=\GK (A_n\t P_m)=2n+m=s,$$
 this proves $(iv)$. $\Box $

A $K$-algebra $A$ is called a {\em central} algebra iff its centre
is $K$.
\begin{corollary}\label{1Wpolc}
\begin{enumerate}
\item The Weyl algebras are the only central  $K$-algebras that
admit a finite set of commuting locally nilpotent derivations with
 generic kernels that have trivial intersection $(=K)$. \item The
polynomial algebras are the only commutative $K$-algebras that
admit a finite set of commuting locally nilpotent derivations with
generic kernels that have trivial intersection $(=K)$.
\end{enumerate}
\end{corollary}

\begin{corollary}\label{eqder}
Let $\d_1, \ldots , \d _s$ be commuting locally nilpotent nonzero
$K$-derivations of an $K$-algebra $A$ satisfying $\cap_{i=1}^s\,
\ker \, \d_i=K$. If $s\neq \GK(A)$ then the kernels of the
derivations $\d_1, \ldots , \d_n$ are not generic.
\end{corollary}

{\it Proof}. Suppose that the kernels of the derivations are
generic. Then by Theorem  \ref{Wpolchar}.(iv), $s=\GK (A)$, a
contradiction. $\Box $

\begin{corollary}\label{1genbassdx=1}
Let $\d_1, \ldots , \d_s$ be a set of commuting locally nilpotent
derivations of an algebra $A$ over a field $K$ of characteristic
zero such that $\d_i (x_j)=\d_{ij}$ for some elements $x_1, \ldots
, x_s$ of $A$. Then $A=C[x_1][x_2; d_2] \cdots [x_s; d_s]$ where
 $C=\cap_{i=1}^s\ker\, \d_i$ and $d_i(x_j):=[x_i,x_j]\in C$.
\end{corollary}

 {\it Proof}. A repeated application of Lemma \ref{dx=1}. $\Box$

The following important  case of this corollary was proved in
\cite{Bass89} (and when $A$ is {\em commutative} the result was
known before [folklore]).

\begin{corollary}\label{bassdx=1}
Let $\d_1, \ldots , \d_s$ be a set of commuting locally nilpotent
derivations of an algebra $A$ over a field $K$ of characteristic
zero such that $\d_i (x_j)=\d_{ij}$ for central elements $x_1,
\ldots , x_s$ of $A$. Then $A=C[x_1, \ldots , x_s]$ is a
polynomial algebra with coefficients from $C=\cap_{i=1}^s\ker\,
\d_i$.
\end{corollary}

\begin{theorem}\label{LN9Nov05}
Let $\s$ be a $K$-algebra endomorphism of the algebra $A:=A_n\t
P_m$ such that $\det (\frac{\der \s (x_{2n+i})}{\der x_{2n+j}})\in
K^*$. Then $\s$ is an algebra automorphism iff the derivations
$\der_1', \ldots , \der_s'$ (see (\ref{dad1}) and (\ref{dad2}))
are locally nilpotent (iff $(\der_i')^{(\deg\,
\s)^{s-1}+1}(x_j)=0$ for all $i,j$, by Theorem \ref{t8Nov05}).
\end{theorem}

{\it Proof}. $(\Rightarrow )$ Obvious.

$(\Leftarrow )$ Since $\der_i'(x_j')=\d_{ij}$ (see (\ref{xsh1}))
for all $i,j=1, \ldots , s$, the derivations $\der_1',\ldots ,
\der_s'$ have generic kernels, then by Corollary
\ref{1genbassdx=1}, $A=\oplus_{\alpha \in
\mathbb{N}^s}C(x')^\alpha$ where $C:=\cap_{i=1}^s\ker (\der_i')$.
It follows from
$$ s=\GK (A)\geq \GK (C)+s$$
that $\GK (C)=0$, which means that every element of the algebra
$C$ is algebraic. Since  scalars are the {\em only} algebraic
elements in $A$ we must have $C=K$ which means  that $A=\s (A)$,
i.e. $\s $ is an automorphism. $\Box $

\begin{proposition}\label{19Nov05}
Let $\s$ be an algebra endomorphism of $A:= A_n\t P_m$ that
satisfies $\det (\frac{\s (x_{2n+i}}{x_{2n+j}})\in K^*$. Then
$N(\der_1', \ldots , \der_s'; A):=\cap_{i=1}^s N(\der_i', A)=\s
(A)$ and $\cap_{i=1}^s\ker_{A}(\der_i')=K$.
\end{proposition}

{\it Proof}. By Corollary \ref{1genbassdx=1}, the intersection is
equal to $N:= C[x_1'][x_2'; d_2']\cdots [ x_s'; d_s']$ where $C:=
\cap_{i=1}^s \ker_A(\der_i')$,  $d_i':= (\ad \, x_i')|_C$, and
$d_i'(x_j'):=[x_i', x_j']\in \{ 0, 1\}$. It follows from
$$ s=\GK (A)\geq \GK (N)\geq \GK (C)+s$$
that $\GK (C)=0$, i.e. each element of $C$ is {\em algebraic} over
$K$, hence $C=K$ (as the only algebraic elements of $A$ are
scalars). $\Box $


\section{Left and right face differential operators}

The algebra $A:=A_n\t P_m$ is {\em self-dual}, i.e. it is
isomorphic to its {\em opposite} algebra $A^{op}$,
$$ A\ra A^{op}, \;\; x_i\mapsto x_{i+n}, \;\; x_{n+i}\mapsto x_i,
\;\; x_{2n+j}\mapsto x_{2n+j}, \;\; i=1, \ldots , n, \;\;
j=1,\ldots , m.$$

\begin{theorem}\label{27Nov05}
Given $\s , \tau \in \Aut_K(A)$ where $A:=A_n\t P_m$ and
$s:=2n+m$. Then
\begin{enumerate}
\item $\s =\tau $ iff $r_i\s =r_i\tau$, $i=1, \ldots , s$, where
$r_i: A\ra A/x_iA$, $a\mapsto a+x_iA$. \item $\s =\tau $ iff
$l_i\s =l_i\tau$, $i=1, \ldots , s$, where $l_i: A\ra A/Ax_i$,
$a\mapsto a+Ax_i$.
\end{enumerate}
\end{theorem}

{\it Proof}. The algebra $A$ is self-dual, so it suffices to prove
that $r_i\s =r_i\tau$, $i=1, \ldots , s$ implies $\s =\tau$.  We
have to show that, for each $k=1, \ldots , s$,
$\s^{-1}(x_k)=\tau^{-1}(x_k)$ (since then $\s^{-1}=\tau^{-1}$
implies $\s =\tau $). $0=r_k (x_k)=r_k\s \s^{-1}(x_k)=r_k\tau
\s^{-1}(x_k)$, hence $\tau \s^{-1}(x_k)\in x_kA$, and so $\tau
\s^{-1}(x_k)=x_ka_k$ for some $a_k\in A$. Applying $\tau^{-1}$ we
obtain $\s^{-1}(x_k)=\tau^{-1}(x_k)b_k$ where
$b_k=\tau^{-1}(a_k)$. By symmetry,
$\tau^{-1}(x_k)=\s^{-1}(x_k)c_k$ for some $c_k\in A$. Now,
$$ \s^{-1}(x_k)=\tau^{-1}(x_k)b_k=\s^{-1}(x_k) c_kb_k \;\; {\rm
and }\;\; \tau^{-1}(x_k)= \s^{-1}(x_k)c_k=\tau^{-1}(x_k) b_kc_k,$$
hence $c_kb_k=b_kc_k=1$ since $A$ is a domain. The only invertible
elements of the algebra $A$ are nonzero scalars, so $c_k, b_k\in
K^*$. Note that the $A/x_iA$ is canonically identified with the
algebra $K\langle x_1, \ldots , \widehat{x}_i, \ldots ,
x_s\rangle$ of type $A_n\t P_{m-1}$ if $x_i$ is central or
otherwise with $A_{n-1}\t P_{m+1}$. Fix $l\neq k$.
$$
x_k=r_l(x_k)=r_l\s\s^{-1}(x_k)=r_l\tau\tau^{-1}(x_k)b_k=b_kx_k,$$
hence $b_k=1$. Therefore, $\s^{-1}(x_k)=\tau^{-1}(x_k)$ for all
$k$, as required.  $\Box $

In the case of the polynomial algebra $A=P_m$, the maps
$r_i=l_i:P_m\ra P_m/(x_i)$ are {\em algebra homomorphisms}, so any
automorphism $\s \in \Aut_K(P_m)$ is uniquely determined by the
algebra epimorphisms $r_i\s :P_m\ra P_m/(x_i)$ {\em which in turn
are uniquely determined by the face polynomials} $\{ r_i\s
(x_j)\}\, | \, i,j=1, \ldots , m\}$ of $\s$ (this is the result of
J. H. McKay and S. S.-S. Wang, \cite{McKayWang88Inv}).

In general situation, $A=A_n\t P_m$, $n\geq 1$, for each $i=1,
\ldots , n$, the maps $r_i$ (resp. $l_i$), are {\em not}  algebra
homomorphisms, they are homomorphisms of {\em right} (resp. {\em
left}) $A$-modules. Note that $A$ is a simple algebra. Note that
$r_{2n+j}=l_{2n+j}$ is an {\em algebra homomorphism} since the
element $x_{2n+j}$ is {\em central} for $j=1, \ldots , m$.

Department of Pure Mathematics

University of Sheffield

Hicks Building

Sheffield S3 7RH

UK

email: v.bavula@sheffield.ac.uk

\end{document}